\documentstyle[12pt]{article}
\newlength{\abstractwidth}
\renewcommand{\thefootnote}{\fnsymbol{footnote}}
\renewcommand{\thanks}[1]{\footnote{#1}} 
\newcommand{\starttext}{
\setcounter{footnote}{0}
\renewcommand{\thefootnote}{\arabic{footnote}}}

\newcommand{\be}{\begin{equation}}
\newcommand{\bea}{\begin{eqnarray}}
\newcommand{\eea}{\end{eqnarray}}
\newcommand{\ee}{\end{equation}}
\newcommand{\N}{{\cal N}}
\newcommand{\<}{\langle}
\renewcommand{\>}{\rangle}
\def\ba{\begin{eqnarray}}
\def\ea{\end{eqnarray}}

\def\N{{\cal N}}

\def\o{\omega}

\def\al{\alpha}
\def\b{\beta}

\def\m{\mu}

\def\o{\omega}

\def\N{\bf N}

\def\i{\infty}

\def\ra{\rightarrow}

\def\N{{\nabla}}

\def\v{\vskip .1in}

\def\[{{\bf [}}
\def\]{{\bf ]}}

\def\pl{\partial}

\begin{document}
\starttext
\baselineskip=18pt
\setcounter{footnote}{0}

\begin{center}
{\Large \bf ON STABILITY AND THE CONVERGENCE OF}\\
\bigskip
{\Large\bf THE K\"AHLER-RICCI FLOW} 
\footnote{Research supported in part by National Science Foundation
grants  DMS-02-45371 and DMS-01-00410}

\bigskip\bigskip

{\large D.H. Phong$^*$ and Jacob Sturm$^{\dagger}$} \\

\bigskip

$^*$ Department of Mathematics\\
Columbia University, New York, NY 10027\\

\v

$^{\dagger}$ Department of Mathematics \\
Rutgers University, Newark, NJ 07102

\end{center}

\baselineskip=15pt
\setcounter{equation}{0}
\setcounter{footnote}{0}

\section{Introduction}
\setcounter{equation}{0}

The normalized K\"ahler-Ricci flow exists for all times,
and converges when the first Chern class is negative or
zero \cite{Ca, Y1}. However, when the first Chern class is positive, there are very few known cases of convergence.
In one complex dimension, Hamilton \cite{H2} used
entropy estimates to show convergence
under the assumption of an initial metric of everywhere positive scalar curvature.
This last assumption was removed later by Chow \cite{Ch}.
In higher dimensions, convergence was established only in the case of positive biholomorphic sectional curvature, first
by X.X. Chen and Tian \cite{CT1, CT2} using Liouville energy functionals, and then by Cao, B.L. Chen, and Zhu \cite{CCZ}
using the recent injectivity radius bound of Perelman \cite{P}.

\medskip
The convergence of the K\"ahler-Ricci flow for a K\"ahler manifold $X$ of positive Chern class can be expected to be a difficult issue, since the limit
would give a K\"ahler-Einstein metric, and not all Fano manifolds
admit such metrics. According to a well-known conjecture of S.T. Yau \cite{Y}, the existence of a K\"ahler-Einstein metric should be equivalent to the stability of $X$ in the sense of geometric invariant theory. It is then an important problem to relate stability to the convergence of the K\"ahler-Ricci flow.

\medskip  
In this paper, we take a first step in this direction. More specifically, we consider the convergence of the K\"ahler-Ricci flow on a compact K\"ahler manifold $(X,\o_0)$ under two assumptions which are both a form of stability:

\medskip
$\ \ $(A) The Mabuchi energy $\nu_{\o_0}(\phi)$ is bounded from below;

\medskip
$\ \ \,$(B) Let $J$ be the complex structure of $X$, viewed as a tensor. Then the $C^\infty$ closure of the orbit of $J$ under the diffeomorphism group of $X$ does not contain any complex structure $J_\infty$ with the property that the
space of holomorphic vector fields
with respect to $J_\infty$ has dimension strictly higher than the dimension of the space of holomorphic vector fields with respect to $J$.

\bigskip

It is well-known by the work of Bando and Mabuchi \cite{BM}
that the existence of a K\"ahler-Einstein metric implies condition (A). This condition is also indirectly related to the notion of $K$-stability \cite{D2, T1}. Indeed, 
if $X$ is imbedded into ${\bf CP}^N$ by sections of the anti-pluricanonical bundle, then $K$-stability can be viewed as
a condition on the asymptotic behavior of the Mabuchi energy
functional $\nu_{\o_0}(\phi)$ along the orbits of $GL(N+1)$.
In particular, it has been shown by Donaldson \cite{D2}
that $K$-stability implies the lower boundedness of
the Mabuchi energy for toric varieties of dimension $n=2$. 

\medskip
The condition (B) is arguably an even more direct manifestation of stability. The stability of a geometric structure should be a
condition insuring that the moduli space of such structures be
Hausdorff. All the tensors $J$ in the same orbit of the diffeomorphism group of $X$ define the same holomorphic structure.
If, as ruled out by condition (B), the dimension of the space of holomorphic vector fields jumps up in the limit, then the limit would be different, and the moduli space of holomorphic structures would not be Hausdorff. 

\bigskip

{\bf Theorem 1} {\it Let $(X,J)$ be a compact complex manifold of dimension $n$. Let $\dot g_{\bar kj}=-R_{\bar kj}+\mu g_{\bar kj}$ be the normalized K\"ahler-Ricci flow, with the initial metric $g_{\bar kj}(0)$ a K\"ahler metric in the first Chern class of $X$. Here $\mu n$ denotes the total scalar curvature.
Assume that the Riemann curvature tensor is uniformly bounded along the flow.

\medskip

{\rm 1}. If condition {\rm (A)} holds,
then we have for any $s\geq 0$
\be
\label{Sobolev}
{\rm lim}_{t\to\infty}||R_{\bar kj}(t)-\mu g_{\bar kj}(t)||_{(s)}\,=\, 0,
\ee
where $||\cdot||_{(s)}$ denotes the Sobolev norm of order $s$ with respect to the metric $g_{\bar kj}(t)$.

\medskip
{\rm 2}. If both conditions {\rm (A)} and {\rm (B)} hold, and if the diameter of $X$ is uniformly bounded along the flow, 
then the K\"ahler-Ricci flow converges exponentially fast in
$C^{\infty}$ to a K\"ahler-Einstein metric.}

\bigskip
Recent as yet unpublished work of Perelman shows that
the scalar curvature and the diameter of the manifold always remain bounded under the K\"ahler-Ricci flow. Thus the diameter and curvature assumptions in Theorem 1 may be significantly less restrictive than they appear at first sight. The following theorem provides a setting where the curvature assumptions in Theorem 1 can be obtained by combining Perelman's result with a relatively mild curvature positivity condition on the initial metric:

\bigskip 

{\bf Theorem 2} {\it Let $X$ be a compact K\"ahler
manifold of dimension $n=2$, which satisfies the conditions 
{\rm (A)} and {\rm (B)}. Assume that
the scalar curvature and the diameter of $X$ are bounded from above along the K\"ahler-Ricci flow with a K\"ahler initial metric in the first Chern class of $X$.
Assume further that

\medskip 

$\ \ $ {\rm (C)} The initial metric $g_{\bar kj}(0)$ has non-negative Ricci curvature and its traceless curvature operator is 2-nonnegative.

\medskip
 
Then the K\"ahler-Ricci flow converges exponentially fast in $C^\infty$ to a K\"ahler-Einstein
metric.}

\bigskip

The notion of $2$-nonnegativity for the Riemann curvature operator was introduced in H. Chen \cite{Chen}.
The condition (C) was introduced in \cite{PS2},
where it was shown in 2 dimensions to be preserved by the K\"ahler-Ricci flow. It should perhaps be mentioned that the case of an initial metric with positive bisectional curvature in arbitrary dimension is another case where the uniform boundedness curvature conditions of Theorem 1 would be satisfied, if we invoke Perelman's unpublished result.
Indeed, by earlier work of Mok \cite{M} and Bando \cite{B},
the positivity of the bisectional curvature is preserved.
Thus the boundedness of the scalar curvature implies the boundedness of the bisectional curvature, and hence of the sectional curvature.

\medskip
Without stability assumptions, the K\"ahler-Ricci flow
is expected to produce either singularities or solitons in
the infinite time limit. For such results, under similar
assumptions of uniform bounds on the curvature, see the recent
works of N. Sesum \cite{Se1,Se2}.

\section{Part 1 of Theorem 1: Sobolev estimates}
\setcounter{equation}{0}

In this section, we prove the first statement in Theorem 1.

\bigskip

$\bullet$ We start with some preliminary considerations about the K\"ahler-Ricci flow and the Mabuchi
functional. Since the initial metric $g_{\bar kj}(0)$ is in $c_1(X)$, and since the K\"ahler-Ricci flow manifestly preserves the first Chern class, we can write at all times
\be
\label{h}
R_{\bar kj}-\mu \,g_{\bar kj}=\pl_j\pl_{\bar k}h
\ee
for some smooth scalar function $h=h(t)$ defined up to a (time-dependent) additive constant. We fix an arbitrary smooth choice of such constants. Since all our estimates ultimately depend only on $\N h$, such choices are immaterial. 
Now $h$ flows according to
\be
\dot h=\Delta h+\mu h+c
\ee
where $c$ is a time-dependent constant, and $\Delta=\N^{\bar k}\N_{\bar k}$. To see this,
it suffices to differentiate the defining equation for $h$
with respect to time. Since $\dot R_{\bar kj}=-\pl_j\pl_{\bar k}
(g^{l\bar p}\dot g_{\bar pl})=\pl_j\pl_{\bar k}R$,
we obtain $\pl_j\pl_{\bar k}R+
\mu\pl_j\pl_{\bar k}h=\pl_j\pl_{\bar k}\dot h$, and hence $R+\mu h+c'=\dot h$,
where $c'$ is a constant. But (\ref{h}) also implies
that $R-\mu n=\Delta h$, so that the desired identity follows
with $c=c'+\mu n$.

\bigskip

$\bullet$ Next, the Mabuchi energy functional $\nu_{\o_0}(\phi)$ is the functional on the space of K\"ahler potentials
$\{\phi;\ g_{\bar kj}=g_{\bar kj}(0)+\pl_j\pl_{\bar k}\phi>0\}$
defined by its variation
\be
\label{mabuchi}
\dot\nu_{\o_0}(\phi)=-{1\over V}\int_X\dot\phi(R-\mu n)\o^n,
\ee
where the K\"ahler form $\o$ is defined by $\o={\sqrt{-1}\over 2\pi}g_{\bar kj}dz^j\wedge d\bar z^k$ and $V=\int_X\o^n$.
To determine $\dot\phi$ in the case of the K\"ahler-Ricci flow,
we rewrite the defining equation $\dot g_{\bar kj}=-R_{\bar kj}+\mu g_{\bar kj}=-\pl_j\pl_{\bar k}h$ in terms of $\phi$. This gives $\dot\phi=-h+c''$, with $c''$ another time-dependent constant. Substituting this in the preceding equation, we get
\be
\dot\nu_\o(\phi)=-{1\over V}\int_X|\nabla h|^2\o^n
\ee

\bigskip

$\bullet$ We return now to the proof of part 1 of Theorem 2 proper.  
Our first step is to show that the lower bound for
the Mabuchi energy functional together with the uniform boundedness of the scalar
curvature implies that
\be
\int_X|\nabla h|^2\o^n\to 0, \ \ \ t\to\infty.
\ee
Now integrating $\dot\nu_\o(\phi)$ and using the lower
bound for $\nu_\o(\phi)$ gives 
\be
\label{below}
{1\over V}\int_0^{T}\,dt\,\int_X|\nabla h|^2\o^n\leq 
\nu_\o(\phi_0)-\nu_\o(\phi_T)
\leq C
\ee
for all $T>0$, which implies that $\int_X|\N h|^2\o^n$ converges to $0$
along some sequence of times tending to $\infty$. To get full
convergence, we consider the flow of $|\nabla h|^2$. 
We have the following Bochner-Kodaira formula 
\be
\Delta|\pl_jh|^2
=g^{j\bar k}\Delta (\pl_jh)\,\overline{\pl_k h}+
g^{j\bar k}\pl_jh\,
\overline{\Delta\pl_kh}
+
R^{j\bar k}\pl_jh\,\overline{\pl_k h}
+
|\bar\nabla\nabla h|^2+|\nabla\nabla h|^2
\ee
Comparing this with the time variation
\be
(|\pl_jh|^2)^{\dot{}}
=
g^{j\bar k}(\pl_jh)^{\dot{}}\,\overline{\pl_kh}
+
g^{j\bar k}\pl_jh\,\overline{(\pl_kh)^{\dot{}}}
+
R^{j\bar k}\pl_jh\,\overline{\pl_kh}
-
\mu|\nabla h|^2
\ee
and the flow $(\pl_jh)^{\dot{}}=\Delta\,(\pl_j h)+\mu\pl_j h$
for $\pl_jh$, we find
\be
\label{flow}
(|\nabla h|^2)^{\dot{}}
-
\Delta|\nabla h|^2
=
-|\bar\nabla\nabla h|^2-|\nabla\nabla h|^2
+
\mu |\nabla h|^2
\ee
Let $Y(t)=\int_X|\nabla h|^2\o^n$. Since $(\o^n)^{\dot{}}=(-R+\mu n)\o^n$, we obtain 
\be
\label{210}
\dot Y
=\mu(n+1)Y-\int_X|\nabla h|^2R\,\o^n
-\int_X|\bar\nabla\nabla h|^2\o^n-
\int_X|\nabla\nabla h|^2\o^n
\ee
Say $|R|\leq C$. Then $\dot Y\leq (\mu(n+1)+C)Y$, and hence
$Y(t)\leq Y(s)e^{(\mu(n+1)+C)(t-s)}$ for all $t\geq s$.
Since the bound (\ref{below}) implies that $\sum_{m=0}^\infty \int_m^{m+1}dt Y(t)<\infty$, there must be a sequence $t_m\in [m,m+1)$ with $Y(t_m)\to 0$. The previous bound implies
$Y(t)\leq Y(t_m)e^{\mu(n+1)+C}$ for all $t\in [m,m+1)$,
and hence $Y(t)\to 0$ as $t\to\infty$.

\bigskip
$\bullet$ The next step is to extend this convergence to the higher derivatives of $h$. Since $Y(t)$ is now known to tend to $0$, the integration of (\ref{210}) gives
\bea
\int_0^\infty dt\int_X|\bar\nabla\nabla h|^2\o^n
+
\int_0^\infty dt\int_X|\nabla\nabla h|^2\o^n
&=&
Y(0)+\mu(n+1)\int_0^\infty dt \int_X|\nabla h|^2\o^n
\nonumber\\
&&
\quad
-
\int_0^\infty dt \int_X|\nabla h|^2R\,\o^n
\eea
This implies that the $L^2$ norms of $\bar\nabla\nabla h$
and $\nabla\nabla h$ tend to $0$ along some subsequence of
times going to infinity. To establish convergence, we need
as previously the flows for $\bar\nabla\nabla h$
and $\nabla\nabla h$. It is convenient to set up a systematic
induction argument as follows. Set 
\bea
&&
h_{\bar K J}=\nabla_{j_s}\cdots\nabla_{j_1}\nabla_{\bar k_r}\cdots\nabla_{\bar k_1}h,
\quad
h_{\bar KJ}\cdot  \overline{h_{\bar KJ}'}
=
g^{L\bar K}g^{J\bar M}
h_{\bar KJ}\overline{h_{\bar LM}'},
\nonumber\\
&&
|\nabla^s\bar\nabla^rh|^2
=
h_{\bar KJ}\cdot  \overline{h_{\bar KJ}},
\quad
g^{L\bar K}g^{J\bar M}=
g^{j_1\bar m_1}\cdots g^{j_s\bar m_s}
g^{l_1\bar k_1}\cdots g^{l_s\bar k_r}
\eea
for $K=(k_1\cdots k_r)$, $J=(j_1\cdots j_s)$. Instead of Bochner-Kodaira formulas for the complex Laplacian $\nabla_p\nabla^p$, it is simpler
to use the real Laplacian $\Delta_{\bf R}=\Delta+\bar\Delta$,
which gives at once
\be
\label{Deltanorm}
{1\over 2}\Delta_{\bf R}|\nabla^s\bar\nabla^rh|^2
=
{1\over 2}\Delta_{\bf R}h_{\bar KJ}\cdot
\overline{h_{\bar KJ}}
+
h_{\bar K J}\overline{{1\over 2}\Delta_{\bf R}h_{\bar KJ}}
+
|\nabla h_{\bar KJ}|^2
+
|\bar\nabla h_{\bar KJ}|^2
\ee
The time evolution of $|\nabla^r\bar\nabla^sh|^2$ is given by
\bea
\label{timenorm}
(|\nabla^s\bar\nabla^rh|^2)^{\dot{}}
&=&
\dot h_{\bar K J}\cdot \overline{h_{\bar K J}}
+
h_{\bar K J}\overline{\dot h_{\bar K J}}
-\mu(r+s)|\nabla^s\bar\nabla^rh|^2\\
&&
\quad
+
\sum_{\al=1}^rR^{l_{\al}\bar k_{\al}}h_{\bar k_1\cdots\bar k_{\alpha}\cdots \bar k_r \,J}
\bar h^{J\bar k_r\cdots}{}_{l_{\alpha}}{}^{\cdots \bar k_1}
+
\sum_{\beta=1}^sR^{j_{\beta}\bar m_{\beta}}
h_{\bar K\,j_1\cdots j_{\beta}\cdots j_s}
\bar h^{j_1\cdots}{}_{\bar m_\beta}{}^{\cdots j_s\bar K}
\nonumber
\eea
Now we need the flow for
the tensor $h_{\bar KJ}$. By induction, we find
\bea
\label{tensorflow}
\dot h_{\bar KJ}
&=&
{1\over 2}\Delta_{\bf R}h_{\bar KJ}+\mu h_{\bar KJ}  
-
{1\over 2}\sum_{\al=1}^rR_{\bar k_{\al}}{}^{\bar m_{\al}}
h_{\bar k_1\cdots\bar m_\al\cdots\bar k_r\,J}
-
{1\over 2}
\sum_{\beta=1}^sR^{m_\beta}{}_{j_{\beta}}h_{\bar Kj_1\cdots m_\beta\cdots j_s}
\nonumber\\
&&
\quad
-
\sum_{1\leq\al<\beta\leq s}R^{m_\al}{}_{j_\al}{}^{m_\beta}{}_{j_\beta}
h_{\bar K j_1\cdots m_{\al}\cdots m_{\beta}\cdots j_s}
-
\sum_{1\leq \al<\beta\leq r}
R_{\bar k_\al}{}^{\bar m_\al}{}_{\bar k_\beta}{}^{\bar m\beta}
h_{\bar k_1\cdots\bar m_\al\cdots\bar m_\beta\cdots \bar k_r\,J}
\nonumber\\
&&
\quad
+
\sum_{\al=1}^r\sum_{\beta=1}^s
R_{\bar k_\al j_\beta}{}^{\bar m_\al n_\beta}
h_{\bar k_1\cdots \bar m_\al\cdots \bar k_r j_1\cdots n_\beta\cdots j_s}
+
\sum_{u=1}^{r+s-2}D^uRiem\,\star\, D^{r+s-u}h
\eea
Here $D$ denotes covariant differentiation in either $j$ or $\bar j$ indices, and $D^uh$ and $D^uRiem$ denote all tensors obtained
by $u$ covariant differentiations of $h$ and of the Riemann
curvature tensor respectively. The $\star$ symbol indicates
general pairings of these tensors. The last term in the above
equation is a lower order term which is actually
absent when $r+s\leq 2$. Assembling the equations
(\ref{Deltanorm}), (\ref{timenorm}), and (\ref{tensorflow}),
we obtain
\bea
(|\nabla^s\bar\nabla^rh|^2)^{\dot{}}
&=&
{1\over 2}\Delta_{\bf R}|\nabla^s\bar\nabla^rh|^2
-|\nabla^{s+1}\bar\nabla^rh|^2
-|\bar\nabla\nabla^s\bar\nabla^rh|^2
+2(\mu-r-s)|\nabla^s\bar\nabla^rh|^2
\nonumber\\
&&
\quad
+2\sum_{\al=1}^r\sum_{\beta=1}^s
R_{\bar k_\al j_\beta}{}^{\bar m_\al n_\beta}
h_{\bar k_1\cdots \bar m_\al\cdots \bar k_r j_1\cdots n_\beta\cdots j_s}\bar h^{\bar KJ}
\nonumber\\
&&
\quad
-
2\sum_{1\leq\al<\beta\leq s}R^{m_\al}{}_{j_\al}{}^{m_\beta}{}_{j_\beta}
h_{\bar K j_1\cdots m_{\al}\cdots m_{\beta}\cdots j_s}
\bar h^{\bar KJ}
\nonumber\\
&&
\quad
-
2\sum_{1\leq \al<\beta\leq r}
R_{\bar k_\al}{}^{\bar m_\al}{}_{\bar k_\beta}{}^{\bar m\beta}
h_{\bar k_1\cdots\bar m_\al\cdots\bar m_\beta\cdots k_r\,J}
\bar h^{\bar KJ}
\nonumber\\
&&
\quad\quad
+
2 
\sum_{u=1}^{r+s-2}D^uRiem\,\star\, D^{r+s-u}h\,
\star\,\overline{\nabla^s\bar\nabla^r h}
\eea
Set $Y_{r,s}(t)=\int_X|\nabla^s\bar\nabla^rh|^2\o^n$. According to \cite{H3, Shi}, the uniform boundedness of the Riemannian curvature
tensor in the Ricci flow implies the uniform boundedness of the
covariant derivatives of the Riemann curvature tensor of any fixed order. The argument applies verbatim to the normalized K\"ahler-Ricci flow. Thus the previous identity implies that
\be
\label{Y}
\dot Y_{r,s}(t)
\leq C_1\,Y_{r,s}(t)+C_2(\int_X|D^{r+s-u}h|^2\o^n)^{1/2}Y_{r,s}^{1/2}(t)
-\int_X|\nabla^{s+1}\bar\nabla^rh|^2\o^n
-\int_X|\bar\nabla\nabla^s\bar\nabla^rh|^2\o^n,
\ee
where a summation over $1\leq u\leq r+s-2$ is understood,
and we have bounded all curvature terms by constants, applied the Cauchy-Schwarz inequality to the lower order terms.
This implies for any $a\geq b$
\bea
\label{Y(a)}
Y_{r,s}(a)-Y_{r,s}(b)
&\leq& C_1\int_b^adt\,Y_{r,s}(t)
+
C_2\int_b^adt\,(\int_X|D^{r+s-u}h|^2\o^n)^{1/2}Y_{r,s}^{1/2}(t)
\nonumber\\
&\leq& C_1\int_b^{\infty}Y_{r,s}(t)
+C_2(\int_b^\infty dt\,\int_X|D^{r+s-u}h|^2\o^n)^{1/2}
(\int_b^\infty dt\,Y_{r,s}(t))^{1/2}\nonumber\\
\eea
We argue now by induction. Assume that
\bea
\int_0^\infty dt\,\int_X|D^vh|^2\o^n &<&\infty,
\quad\quad
{\rm for}\ v\leq r+s
\nonumber\\
\int_X|D^vh|^2\o^n&\rightarrow& 0,
\quad\quad
{\rm for}\ v<r+s.
\eea
Then $Y_{r,s}(t)$ is in particular integrable on $[0,\infty)$, and arguing as before, we can choose $b_m\in[m,m+1)$ with $Y_{r,s}(b_m)\to 0$. The estimate (\ref{Y(a)}) applied
with $b=b_m$ and $m$ large enough shows that $Y_{r,s}(t)\to 0$.
Since any covariant derivative of $h$ of order $r+s$ differs from covariant derivatives of the form $\nabla^{s}\bar\nabla^rh$ and $\nabla^s\bar\nabla^rh$
by $D^uRiem\star D^{r+s-u}h$ with $u\geq 1$, and since all derivatives of the Riemann curvature tensor are bounded,
it follows from the second induction hypothesis and the
fact that $Y_{r,s}(t)\to 0$ that $\int_X|D^{r+s}h|^2\o^n\to 0$.

\medskip
To establish the first induction hypothesis at order $r+s+1$,
we return to the equation (\ref{Y}). Integrating from
$0$ to $\infty$, and applying again the Cauchy-Schwarz inequality, we obtain
\bea
&&
\int_0^\infty dt\,\int_X|\nabla^{s+1}\bar\nabla^rh|^2\o^n
+
\int_0^\infty\,\int_X|\bar\nabla\nabla^s\bar\nabla^rh|^2\o^n
\\
&&
\quad
\leq
Y_{r,s}(0)+C_1\int_0^\infty dt\,Y_{r,s}(t)
+
C_2
(\int_0^\infty dt\,\int_X|D^{r+s-u}h|^2\o^n)^{1/2}
(\int_0^\infty dt\,Y_{r,s}(t))^{1/2}
\nonumber
\eea 
This implies that the $L^2$ norms of $\nabla^{s+1}\bar\nabla^rh$
and $\bar\nabla\nabla^s\bar\nabla^rh$ are integrable with respect to time on $[0,\infty)$. Using the first induction hypothesis and again the uniform boundedness of the Riemann curvature tensor and its derivatives, we can deduce that $\int_X|D^{r+s+1}h|^2\o^n$
is integrable with respect to time on $[0,\infty)$.
This completes the induction argument. Since $R_{\bar kj}-\mu g_{\bar ki}$ is given by $\pl_{\bar k}\pl_jh$,
the $L^2$ convergence to $0$ of all covariant derivatives of $h$
implies the convergence to $0$ of all Sobolev norms of $R_{\bar kj}-\mu g_{\bar kj}$. The proof of the first
part of Theorem 1 is complete.

\section{Part 2 of Theorem 1: convergence of the metrics}
\setcounter{equation}{0}

Assuming now condition (B) and the uniform boundedness of the diameter of $X$, we wish to establish the convergence of the metrics $g_{\bar kj}(t)$ as $t\to\infty$. The uniform boundedness
of the curvature together with the uniform boundedness of the diameter imply the uniform boundedness from below of the injectivity radius \cite{C1}. Since the volume is fixed under the normalized flow, the uniform control of the diameter and the injectivity radius imply the uniform control of the Sobolev constant. Thus the equation (\ref{Sobolev}) now implies
\be
\label{C0}
{\rm sup}_X|D^p\dot g_{\bar kj}(t)|_t
={\rm sup}_X|D^p(R_{\bar kj}(t)-\mu g_{\bar kj}(t))|_t\to 0,
\ee
for any fixed integer $p$,
where we have introduced the lower index $t$ to stress that the
norms are taken with respect to the metric $g_{\bar kj}(t)$.
However, the convergence to $0$ of $|\dot g_{\bar kj}(t)|_t$
does not guarantee the convergence of the metrics themselves.
In fact, it does not even guarantee that they are uniformly
equivalent. 

\v
We also note that the uniform boundedness of the volume, diameter, injectivity radius, and curvature implies that
the metrics $g_{\bar kj}(t)$ have uniform bounded geometry,
in the sense of Gromov \cite{G}, in fact uniform bounded $C^\infty$ geometry, since all covariant derivatives of the curvature are also bounded uniformly. Thus, by passing to a subsequence and applying the $C^\infty$ version of Gromov compactness due to Hamilton \cite{H4}, we can find diffeomorphisms $F_{t_j}$ so that
the pull-backs $(F_{t_j})_{*}(g(t_j))$ converge in $C^\infty$.
But we have no control over the diffeomorphisms $F_{t_j}$
and cannot deduce from this the convergence of the metrics $g_{\bar kj}(t_j)$ themselves. As we saw earlier, this issue of diffeomomorphisms underlies the notion of stability, and it
appears now central to the problem of convergence of the flow.
We shall see later, however, that Gromov compactness can be
put to good use in the proper context.

\bigskip
$\bullet$ To overcome these difficulties, we shall establish the exponential decay of $|\dot g_{\bar kj}(t)|_t$.
Let $Y=\int_X|\nabla h|^2\o^n$ as before. Then we have
\be
\label{key}
\dot Y
=-\int_X|\nabla h|^2(R-\mu n)\o^n-\int_X\N^jh\overline{\N^kh}(R_{\bar kj}-\mu g_{\bar kj})\o^n 
-2\int_X|\bar\nabla \bar\nabla h|^2\o^n
\ee
This follows from the equation (\ref{flow}),
the fact that $(\o^n)^{\dot{}}=(R-\mu n)\o^n$,
and the Bochner-Kodaira formula for vector fields $V=(V^j)$,
\be
||\nabla V||^2=||\bar\nabla V||^2+\int_X R_{\bar kj}V^j\bar V^{\bar k}\,\o^n
\ee
applied to $V^j=g^{j\bar k}\pl_{\bar k}h=\N^j h$
(note that $\N^jh=(\bar\N h)^j)$.
The identity (\ref{key}) shows that the only chance of getting exponential decay is by obtaining a strictly positive lower bound for $\int_X|\bar\nabla \bar\nabla h|^2\o^n=\int_X|\bar\partial V|^2\o^n$.

\bigskip
$\bullet$
Let $\lambda_t$ be the lowest strictly positive eigenvalue
of the complex Laplacian $\Delta_t=g^{j\bar k}\N^j\N^{\bar k}
=-\bar\pl^{\dagger}\bar\pl$ on $T^{1,0}(X)$ vector fields. 
By the elliptic theory,
we have
\be
\lambda_t\,||V-\pi_t V||^2\leq \int_X|\bar\nabla V|^2\o^n,
\ee
where $\pi_t$ is the orthogonal projection with respect to $g_{\bar kj}(t)$ on the space $H^0(X,T^{1,0})$ of holomorphic vector fields. When $V$ is the gradient of the function $h$ as in our present case, it turns out that $||\pi_tV||^2$ is exactly the
Futaki invariant of the manifold $X$, applied to $\pi_tV$.
Recall that the Futaki invariant $Fut$ is the character defined on the space of holomorphic vector fields by \cite{F}
\be
\label{futaki}
Fut(W)=\int_X (Wh)\o^n,
\quad\quad W\in H^0(X,T^{1,0}).
\ee
It is independent of the choice of metric within the K\"ahler class. The main observation for our purposes is
\be
||\pi_t V||^2=Fut(\pi_t V)
\ee
when $V^j=\nabla^j h$. To see this, we note $\<\pi_t V,\pi_t V\>
=\<V,\pi_t^2V\>=\<V,\pi_t V\>=\<\bar\N h,\pi_t V\>=(\pi_t V)(h)
=Fut(\pi_t V)$. Thus the inequality (\ref{key}) implies
\be
\dot Y\leq -2\,\lambda_t\,Y+2\,\lambda_t\,Fut(\pi_t(\N^j h))
-\int_X|\nabla h|^2(R-\mu n)\o^n-\int_X\N^jh\,\overline{\N^kh}(R_{\bar kj}-\mu g_{\bar kj})\o^n.
\ee
This key inequality holds in all generality for the normalized K\"ahler-Ricci flow.

\bigskip
$\bullet$ Assume now all the conditions stated for part 2 of Theorem 1. Then the lower boundedness of the Mabuchi energy functional implies that the Futaki invariant vanishes identically
(see \cite{T2}). Furthermore, $|\N^jh\,\overline{\N^kh}(R_{\bar kj}-\mu g_{\bar kj})|\leq |\N h|^2|R_{\bar kj}-\mu g_{\bar kj}|_t$, so that the convergence of $|R_{\bar kj}-\mu g_{\bar kj}|_t$ to $0$ established in (\ref{Sobolev}) implies
that for any $\epsilon>0$, we have
\be
\dot Y\leq (-2\,\lambda_t+\epsilon)\,Y,
\quad\quad t\in [T_\epsilon,\infty),
\ee
for $T_\epsilon$ large enough.
Thus establishing exponential convergence reduces to showing that $\lambda_t$ is uniformly bounded from below by a positive constant.

\bigskip
$\bullet$ The bounds from below for $\lambda_t$ do not appear to be accessible by Bochner-Kodaira techniques, as these
apply to a negative bundle instead of a positive one, as is here the case. Rather, they should reflect the bounded geometry, just as in the case of the lowest eigenvalue for the scalar Laplacian.
In that case, we can either construct explicitly the Green's function as in \cite{A} or make use of estimates for the lowest eigenvalue such as Cheeger's \cite{C} in terms of the isoperimetric and Sobolev constants \cite{CL}. However, no simple characterization of $\lambda_t$ seems available in the case of vector fields. Instead, we shall establish the following estimate
by using a complex version of the Gromov compactness theorem:

\bigskip
{\bf Theorem 3}. {\it Let $X$ be a compact, complex manifold
of dimension $n$. Assume that its complex structure $J$ is stable
in the sense that it satisfies the condition (B) stated in the
Introduction. Fix $V>0$, $D>0$, $\delta>0$, and constants $C_k$.
Then there exists an integer $N$ and a
constant $C=C(V,D,\delta,C_k,n,N)>0$ 
such that
\be
C\,||W||^2\leq ||\bar\pl W||^2,
\quad\quad
W\perp H^0(X,T^{1,0}),
\ee
for all K\"ahler metrics $g$ on $X$ whose volumes and diameters are bounded above by $V$ and $D$ respectively, whose injectivity radius is bounded from below by $\delta>0$, and the $k$-th derivatives of whose curvature tensors are uniformly bounded by $C_k$, for all $k\leq N$.}

\bigskip
Assuming Theorem 3 for the moment, we deduce that there exists a positive constant $c$ so that $\lambda_t>2c$ for all $t$. Thus for $t$ large enough, $Y(t)$
satisfies the differential inequality
\be
\dot Y(t)\leq -c\,Y(t)
\ee
from which it follows that $Y(t)$ decreases exponentially fast
\be
\int_X|\nabla h|^2\o^n\leq C\,e^{-c\,t}.
\ee

\bigskip
$\bullet$ Once the exponential decay of the $L^2$ norm of $\N h$,
has been established, it is not difficult to deduce
the exponential decay of the $L^2$ norms $Y_{r,s}$ of $\bar\N^r\N^s h$,
where all norms are taken with respect to the metric $g_{\bar kj}(t)$. For example, the inequality (\ref{Y}) implies that
\be
\dot Y_{r,s}(t)
\leq (C_1+{1\over 2})\,Y_{r,s}(t)+{1\over 2}\int_X|D^{r+s-u}h|^2\o^n
-\int_X|\nabla^{s+1}\bar\nabla^rh|^2\o^n
-\int_X|\bar\nabla\nabla^s\bar\nabla^rh|^2\o^n,
\ee  
and on the right hand side, a summation over all indices
$u$ in the range $1\leq u\leq r+s-2$ is understood.
For any
$\epsilon>0$, there exists
a constant $c(\epsilon)$ independent of $t$ so that
\be
Y_{r,s}(t)\leq \epsilon 
\big(\int_X|\nabla^{s+1}\bar\nabla^rh|^2\o^n
+\int_X|\bar\nabla\nabla^s\bar\nabla^rh|^2\o^n
\big)
+
c(\epsilon)\,
\int |D^{r+s-1}h|^2\o^n.
\ee
By choosing $\epsilon$ small enough, we can deduce that
\be
\dot Y_{r,s}(t)
\leq
\sum_{1\leq u\leq r+s-1}\int_X|D^{r+s-u}h|^2\o^n
\leq C\ e^{-c\,t},
\ee
where in the last inequality, we have assumed by induction that all $L^2$ norms of $D^{r+s-u}h$ decay exponentially. Since we already know that $Y_{r,s}(t)\to 0$ as $t\to\infty$, we may integrate between $t$ and $\infty$, and see that $Y_{r,s}(t)$
decays exponentially also.

\bigskip
$\bullet$ From the exponential decay of the $L^2$ norms of
$\N^p h$ for all $p\geq 1$, we deduce from the Sobolev imbedding theorem with uniform constants that we have exponential
decay of the $C^k$ norms
\be
{\rm sup}_X |\N^kh|_t^2\leq C_k \ e^{-c\,t}.
\ee

\bigskip
$\bullet$ The next step is to show that all metrics $g_{\bar kj}(t)$ are uniformly equivalent, that is, bounded by one another up to a constant independent of $t$. 
According to a lemma of Hamilton (\cite{H1}, Lemma 14.2),
it suffices to show that
\be
\int_T^\infty {\rm sup}_X|\dot g_{\bar kj}|_t\,dt<\infty.
\ee
Since $\dot g_{\bar kj}(t)=\pl_{\bar k}\pl_j h$, the preceding result implies that
$|\dot g_{\bar kj}|_t$ decays exponentially as $t\to \infty$.
This implies the desired inequality, and hence all metrics
$g_{\bar kj}(t)$ are uniformly equivalent.

\bigskip
$\bullet$ Since all metrics $g_{\bar kj}(t)$ are now known to be equivalent, we can now write for some strictly positive constant $c$ 
\bea
|g_{\bar kj}(T)W^j\bar W^{\bar k}
-
g_{\bar kj}(S)W^j\bar W^{\bar k}|
&\leq&
\int_S^T{\rm sup}_X|\dot g_{\bar kj}|_t|W|_t^2\ dt
\leq C\, |W|_{t=0}^2\int_S^T{\rm sup}_X|\dot g_{\bar kj}|_t\ dt,
\nonumber\\
&\leq& C'\,|W|_{t=0}^2\ (e^{-c\,T}-e^{-c\,S}).
\eea 
This tends to $0$ exponentially as $S,T\to\infty$. Thus the metrics $g_{\bar kj}(t)$ converge exponentially fast as $t\to\infty$ to some metric
$g_{\bar kj}(\infty)$, which is also equivalent to all the $g_{\bar kj}(t)$'s. Iterating the arguments shows that the convergence is in $C^\infty$. Since $\pl_{\bar k}\pl_jh$
tends to $0$, the metric $g_{\bar kj}(\infty)$
is clearly K\"ahler-Einstein.
The proof of Theorem 1 is complete. Q.E.D.

\section{Lower bounds for the $\bar\pl$ operator on vector fields and Gromov compactness}
\setcounter{equation}{0}

It remains to prove Theorem 3. This theorem is essentially a consequence of the following K\"ahler version of the Gromov compactness theorem, combined with some elementary perturbation theory for Laplacians:

\bigskip
{\bf Theorem 4.} {\it Let $X$ be a compact smooth manifold.
Let $(g(t),J(t))$ be any sequence of metrics $g(t)$ and complex structures $J(t)$ on $X$ such that $g(t)$ is K\"ahler with respect to $J(t)$. Assume that the $g(t)$'s
have bounded geometry, in the sense that their volumes, diameters, curvatures, and covariant derivatives of their curvature tensor are all bounded from above, and their injectivity radii are all bounded from below. Then there exists 
a subsequence $t_j$, and a sequence of diffeomorphisms $F_j:X\to X$ such that the pull-back metrics $\tilde g(t_j)= F_j^*g(t_j)$ converge in $C^\infty$ to a smooth metric $g(\infty)$,
the pull-back complex structure tensors $\tilde J(t_j)=F_j^*J(t_j)$ converge in $C^\infty$ to an integrable complex structure tensor
$\tilde J(\infty)$. Furthermore, the metric $\tilde g(\infty)$ is
K\"ahler with respect to the complex structure $\tilde J(\infty)$.}

\bigskip
{\it Proof of Theorem 4}: The $C^\infty$ part of the theorem, without reference to complex structures and K\"ahler forms, is actually the version of the Gromov compactness theorem established by Hamilton \cite{H3}, where $C^k$ uniform bounds
on the curvature are assumed, instead of just $C^0$ bounds as in
the original version of Gromov, Peters \cite{P}, Greene and Wu \cite{GW} (Hamilton's version is even more difficult, because no diameter bound is assumed. In that case, the diffeomorphisms $F_{t_j}$ map a sequence of exhausting compact subsets of $X$
into a sequence of exhausting compact subsets of a limiting manifold $\tilde X$, which may not be compact. One also needs to
choose reference points $P_t$ and reference frames at these points.) Thus, passing to a subsequence, we assume that there exist
diffeomorphisms so that $\tilde g(t)$ converges, and concentrate
on finding a subsequence $t_j$ so that the complex structures $\tilde J(t_j)$ converges also.

\v
First we recall the definitions: To say that a metric $g$ is compatible with a complex structure $J$
means that $g(u,v)=g(Ju,Jv)$
for all $u,v\in TM$. In local coordinates,
$g_{ij}u^iv^j\ = \ g_{kl}J_i^ku^iJ_j^lv^i$,
in other words
\be
\label{hermitian}
g_{ij}\ = \ J_i^kg_{kl}J_j^l
\ee
To say that $J$ is a complex structure is to say that $J^2=-I$
and its Nijenhuis tensor vanishes.
\v

Let $g=g_\i$ and let $\N$ be the Riemannian connection associated to $g$. It suffices to show that there are constants $C_\al$ such that
\be
\label{2}
||\N^\al J_t||_g \ \leq \ C_\al \ \ \hbox{for all $t$}.
\ee
We do this first in the case $\al=0$: Since $g_t\ra g$, it suffices to prove
\be
\label{3}
|| J_t||_{g_t} \ \leq \ C_\al \ \ \hbox{for all $t$}
\ee
Working in normal coordinates for $g_\m$, the equation (\ref{hermitian}) implies, for each
$i$, that
$$ 1= \sum_k J_i^kJ^k_i
$$
Thus $||J_t||_{g_\m}=n$, and this proves (\ref{3}).
\v
Now we prove (\ref{2}) by induction: It's true when $\al=0$.
Since $\N_t J_t=0$ (this is the definition of ``K\"ahler") we have
\be
\label{difference}
\N^\al J_t\ = \ \N^{\al-1}(\N-\N_t)J_t
\ee
Let $H_t=\N-\N_t$. Then $(H_t)g_t=\N g_t$. In other words,
\be
(H_t)^p_{ki}(g_t)_{pj}+(H_t)^p_{kj}(g_t)_{pi}\ = \ \N_k(g_t)_{ij}
\ee
Permuting the indices gives
\bea
(H_t)^p_{ij}(g_t)_{pk}+(H_t)^p_{ik}(g_t)_{pj}\ &=& \ \N_i(g_t)_{jk}
\nonumber\\
(H_t)^p_{jk}(g_t)_{pi}+(H_t)^p_{ji}(g_t)_{pk}\ &=& \ \N_j(g_t)_{ki}
\eea
Thus
\be
2(H_t)^p_{ij}\ = \
(g_t)^{pk}[\N_j(g_t)_{ki}+\N_i(g_t)_{jk}-\N_k(g_t)_{ij}]
\ee
This shows that $H_t$ and all its derivatives are uniformly bounded:
indeed,  $H_t$ converges in $C^\i$ to $H_\i=0$. If follows from (\ref{difference}) that
$J_t$ and its derivatives are bounded. Thus a subsequence converges, and the limit $J_\infty$ is clearly a complex structure. Since
$H_t J_t=\N J_t $, we get $H_\i J_\i=\N J_\i$, but $H_\i=0$. 
Thus $J_\i$ is K\"ahler, and
the proof of Theorem 4 is complete. Q.E.D.

\bigskip
{\it Proof of Theorem 3}: The main step in the proof is to show that, if $(g(t),J(t))$ are a sequence of K\"ahler metrics which converges in $C^\infty$ to $(g(\infty),J(\infty))$, where $g(\infty)$ is a K\"ahler metric with respect to the complex structure $J(\infty)$, and if the dimension of the space of holomorphic vector fields is the same for all $N\leq t\leq\infty$, then
\be
\label{eigenvaluelimit}
{\rm lim}_{t\to\infty}\lambda_t=\lambda_\infty,
\ee
where $\lambda_t$ is the lowest strictly positive eigenvalue
of the Laplacian $\Delta_t=\N_t^j\N_{\bar j}$ on the space
$T_t^{1,0}(X)$ of complex tangent vectors with respect to
$(g(t),J(t))$. 
Assuming this for the moment, Theorem 3 can be proven by contradiction: if it does not hold, then there exists a sequence of metrics $g(t)$ with $\lambda_t\to 0$.
Passing to a subsequence, we can apply 
Theorem 4, and obtain diffeomorphisms $F_t$
so that, the metrics 
$\tilde g_t=(F_t)_*(g_t)$ and complex structures $\tilde J(t)=(F_t)_*(J)$ converge in $C^\infty$ to a metric
$\tilde g_\infty$ and complex structure $\tilde J(\infty)$. By the preceding inequality,
the lowest eigenvalues of $\tilde g(t)$
tend then to a strictly positive limit. But 
$F_t$ is a biholomorphic isometry between the space $X$ equipped with the K\"ahler structure $(g(t),J(t)$ and the K\"ahler structure $(\tilde g(t),\tilde J(t))$. Thus the eigenvalues
of $(\tilde g(t),\tilde J(t))$ are the same as the eigenvalues of $(g(t),J(t))$. This contradiction proves Theorem 3.

\bigskip

$\bullet$ We turn now to the proof of the eigenvalue limit
(\ref{eigenvaluelimit}). Let $||\cdot||_{H_t^{(s)}}$ be the Sobolev norm of order $s$ on $E$, taken with respect to the metric $g_t$. Since the metrics $g(t)$ converge, we have the following inequalities
\bea
&&|\<U,V\>_t-\<U,V\>_\infty|
\leq c_t\,||U||_\infty ||V||_\infty,
\quad\quad
U,V\in C^\infty(X,T(X)),
\nonumber\\
&&
C_t^{-1}||V||_{H_t^{(s)}}
\leq
||V||_{H_\infty^{(s)}}
\leq
C_t
||V||_{H_t^{(s)}},
\quad\quad V\in C^\infty(X,T(X)).
\eea
with constants $c_t\to 0$, $C_t\to 1$ as $t\to\infty$.
Here $\<\cdot,\cdot\>_t$ denotes the inner product with respect to $g_t$. Furthermore, there exists constants $C$ independent of $t$ so that the elliptic a priori estimate for $\Delta_t$ holds uniformly in $t$:
\be
||V||_{H_t^{(1)}}
\leq
C(\<\Delta_tV,V\>_t+||\cdot||_{H_t^{(0)}}^2),
\quad\quad V\in C^{\infty}(X,T^{1,0}(X)).
\ee
Let $\{\phi_t^{(\al)}\}_{1\leq\al\leq N}$ be an orthonormal set of eigenvectors for $\Delta_t$
\be
\Delta_t\phi_t^{(\al)}=0,
\quad\quad 
\<\phi_t^{(\al)},\phi_t^{(\beta)}\>=\delta^{\al\b}.
\ee
For each $\al$, the uniform a priori estimate implies that $\phi_t^{(\al)}$ is uniformly bounded in $H_t^{(1)}$.
But the uniform equivalence of all Sobolev norms $||\cdot||_t^{(1)}$ to $||\cdot||_\infty^{(1)}$ implies that
$\phi_t^{(\al)}$ is uniformly bounded in $H_\infty^{(1)}$.
By Rellich's lemma, it follows that there is a subsequence $\phi_{t_j}^{(\al)}$ which converges in $L_\infty^2(X)$.
This implies
$\<\phi_\infty^{(\al)},\phi_\infty^{(\beta)}\>=
lim_{j\to\infty}
\<\phi_t^{(\al)},\phi_t^{(\beta)}\>_t
=
\delta^{\al\b}$.
Furthermore, $\Delta_\infty\phi_\infty^{(\al)}=
lim_{j\to\infty}\Delta_t\phi_t^{(\al)}=0$
in the sense of distributions, which implies by elliptic regularity that
$\phi_\infty^{(\al)}$ is actually smooth, and the
equation holds in the standard sense. Thus we have shown that,
there is a {\it subsequence} $t_j\to\infty$ so that
for each $\al$, the sequence $\phi_{t_j}^{\al}$ converges 
in $L^2$ to an orthonormal set 
$\phi_\infty^{(\al)}\in H^0_\infty(X,T^{1,0}(X))$.

\bigskip
$\bullet$ We make use now of the assumption that the dimensions
of $K_t=Ker(\Delta_t)$ and $K_\infty=Ker(\Delta_\infty)$ are the same. Thus, by passing to a subsequence,
we can assume the existence of 
orthonormal bases $\{\phi_t^{(\al)}\}$
for $K_t=Ker(\Delta_t)$ for $1\leq t\leq\infty$, with
$\phi_t^{(\al)}$ converging to $\phi_\infty^{(\al)}$ as $t\to\infty$. 
We now show that
\be
{\rm liminf}_{t\to\infty}\lambda_t\geq \,\lambda_\infty
\ee
Let $K_t^{\perp}$ be the orthogonal complement of
$K_t$ in $L^2(X,T_t^{1,0}(X))$, with respect to the metric $g_t$.
Let $\psi_t\in K_t^\perp$ be a lowest eigenfunction of $\Delta_t$,
\be
\Delta_t\psi_t=\lambda_t\psi_t,
\quad\quad\psi_t\in K_t^\perp,
\quad
||\psi_t||_t^2=1.
\ee
Fix any $\epsilon>0$. Assume that there exists a sequence $t_j\to\infty$ so that
\be
\label{liminf}
\lambda_t\leq (1-\epsilon)\lambda_\infty
\ee
for all $t_j$. We abbreviate $t_j$ by $t$ for the sake of notational simplicity. Then $||\Delta_t\psi_t||_{L_t^2}
=\lambda_t$ is bounded, and the uniform elliptic
a priori estimate implies that $||\psi_t||_{H_t^{(2)}}$
is uniformly bounded. In fact, the same argument applied to
$\Delta_t^2$ and its corresponding uniform elliptic a priori
estimate implies that $||\psi_t||_{H_t^{(4)}}$ is uniformly
bounded. Thus we may assume that $\psi_t$ converges in $H_\infty^{(2)}$. Clearly the limit is in $K_\infty^\perp$.
Now let $\Pi$ denote the orthogonal projection from
$T(X)$ to $K_\infty^\perp$. It follows immediately from the
convergence of $(g(t),J(t))$ to $(g(\infty),J(\infty))$,
and the convergence of the orthonormal basis $\{\phi_t^{(\al)}\}$
of $K_t$ to the orthonormal basis $\{\phi_\infty^{(\al)}\}$
of $K_\infty$ that
\be
||\Pi \psi_t-\psi_t||_{H_\infty^{(2)}}\to 0,
\ee
for any $\psi_t\in K_t^\perp$ converging to a vector field in $K_\infty^\perp$. 
Since $\Delta_t-\Delta_\infty$ is a differential operator of second order whose coefficients tend to $0$ in $C^\infty$, we have $||\Delta_t-\Delta_\infty||_{Hom(H_{\infty}^{(2)},
H_\infty^{(0)})}\to 0$, where the norm denotes the operator norm
from $H_\infty^{(2)}$ to $H_\infty^{(0)}$.
Thus
\be
\<\Delta_t\psi_t,\psi_t\>_t
=
\<\Delta_\infty(\Pi\psi_t),\Pi\psi_t\>_\infty-o(1)
\geq
\lambda_\infty\,||\Pi\psi_t||_\infty^2-o(1),
\ee
where $o(1)$ denotes terms tending to $0$ as $t\to\infty$.
But $||\Pi\psi_t||_\infty^2\to 1$ as $t\to\infty$. This contradicts (\ref{liminf}) for $t$ large enough.
Our statement about the
${\rm liminf}\,\lambda_t$ follows at once.

\bigskip

$\bullet$ Similarly, we can show that
${\rm limsup}\,\lambda_t\leq\lambda_\infty$,
even without the assumption about the ranks of $K_t$ not jumping up in the limit. Indeed, for any fixed $\epsilon>0$, choose $\psi\in Ker(\Delta_\infty)^\perp$ with
$||\psi||_\infty=1$ and
$\lambda_\infty\leq \<\Delta_\infty\psi,\psi\>_\infty
\leq \lambda_\infty+\epsilon$.
If $\Pi_t$ denotes the orthogonal projection from $T(X)$ to
$K_t^\perp$, then it is easy to see that $\Pi_t\psi\in K_t^\perp$
satisfies $\<\Delta_t(\Pi_t\psi),\Pi_t\psi\>_t
\to\<\Delta_\infty\psi,\psi\>_\infty$.
Since $\lambda_t\leq \<\Delta_t\psi_t,\psi_t\>_t$,
this gives the desired estimate. The proof of Theorem 3 is complete. Q.E.D.

\section{Proof of Theorem 2}
\setcounter{equation}{0}

Under the assumption of non-negativity of the Ricci curvature and 2-nonnegativity of the traceless curvature operator for the initial metric, the non-negativity of 
the Ricci curvature is preserved for all times
\cite{PS2}. Thus the boundedness of the
scalar curvature implies the boundedness of the eigenvalues
and hence of the Ricci curvature.

\medskip
Next, we consider the curvature tensor $R_{\bar ab\bar cd}$.
On a K\"ahler manifold, the Riemann curvature operator can be viewed as an operator on the space of real $(1,1)$-forms. The
traceless curvature operator $Op(S)$ is the projection of the Riemann curvature operator on the subspace of
traceless real (1,1)-forms (see \cite{PS2}, eq.(2.8)). Now, under the assumption (a),
the eigenvalues $m_1\leq m_2\leq m_3$ of the traceless 
curvature operator $Op(S)$
are bounded: indeed, if they are all non-negative, then this follows
immediately from the boundedness of the scalar curvature since
$m_1+m_2+m_3=R/2$.
Otherwise, the 2-nonnegativity implies that at most one eigenvalue
$m_1$ is negative, and that $0\leq m_1+m_2$. But then
\be
{1\over 2}R=(m_1+m_2)+m_3\geq m_3\geq 0,
\ee
and hence $m_3$ is uniformly bounded. Since $0\leq m_2\leq m_3$, so is $m_2$.
Finally,
\be
0\leq {1\over 2}R=m_1+m_2+m_3
\Rightarrow
|m_1|\leq m_2+m_3
\ee
and thus $|m_1|$ is also uniformly bounded. 
The boundedness of the eigenvalues of $Op(S)$ implies
that its entries are also uniformly bounded. This is because the matrix for $Op(S)$ is symmetric, and thus diagonalizable by
unitary matrices. Next, the curvature operator $R_{\bar ab\bar cd}=Op(R)$ can be written as
\be
Op(R)=\pmatrix{R/2 & S\cr S^t& Op(S)\cr}
\ee
where $S$ is the traceless part of the Ricci curvature
$S_{\bar ab}=R_{\bar ab}-R\delta_{\bar ab}/n$.
We deduce that the entries of $Op(R)$
are all bounded. Hence its eigenvalues are also bounded, establishing the uniform boundedness of the Riemann
curvature tensor $R_{\bar ab\bar cd}$.
Theorem 2 follows now from Theorem 1. Q.E.D.

\section{Remarks}
\setcounter{equation}{0}

$\bullet$ It is clear from the proof of Theorem 1 that, for any $\epsilon>0$, the rate of exponential convergence can be taken to be $\lambda_{\infty}-\epsilon$, where $\lambda_{\infty}$ is the lowest
strictly positive eigenvalue of the complex Laplacian
on the space $T^{1,0}(X)$ of vector fields with respect to
the K\"ahler-Einstein metric $g_{\bar kj}(\infty)$. 

\bigskip

$\bullet$ Several specific notions of stability have been by now proposed in the literature \cite{D1, D2, PT, PS1, RT, T1}. 
At the present time, the relations between these various notions are still obscure. Nor has any precise relation
between any of them and the convergence of the K\"ahler-Ricci flow been as yet proved.
This is clearly an important direction for further investigation,
and it can be hoped that the methods of the present paper
would be useful.

\bigskip 

$\bullet$ From the proof of part 1 of Theorem 1, it is clear that
the assumption that $R\geq 0$ (but not necessarily uniformly
bounded from above), combined with lower bounds for the Mabuchi energy also suffices to show that $\int_X|Dh|^2\o^n\to 0$ and $\int_0^\infty dt\,\int_X|D^2 h|^2\o^n<\infty$.

\bigskip 

$\bullet$ The preceding inequality implies that           $\int_0^\infty\,\int_X(R-\mu n)^2\o^n<\infty$. This inequality
was instrumental in \cite{CT1}, where it was
established under the stronger assumption
of positive biholomorphic sectional curvature,
using the Moser-Trudinger inequality and Liouville
functionals.

\bigskip

$\bullet$ The flow for $|D^rh|^2$ for general $r$ requires
bounds for the full Riemann curvature tensor and its derivatives. However, certain lower order derivatives can still be bounded
under the assumption of lower bounds for the Mabuchi functional
and weaker curvature assumptions. For example, under the weaker assumption that $|R|$ remains bounded, we still have
\be
\int_X|\nabla\bar\nabla h|^2\o^n\to 0.
\ee
Indeed, an integration by parts shows that the $L^2$ norm of $\nabla\bar\nabla h$ is the same as the $L^2$ norm of $\Delta h$.
The flow of $\Delta h$ is easily derived from that of
$h_{\bar kj}$, which is given in (\ref{tensorflow}),
with no lower order terms as we had noted earlier. We find
$(\Delta h)^{\dot{}}=\Delta(\Delta h)+R^{j\bar k}h_{\bar kj}$,
and hence
\be
(|\Delta h|^2)^{\dot{}}
=
\Delta(|\Delta h|^2)-2|\nabla\Delta h|^2+2|\nabla\bar\nabla h|^2\Delta h
+2\mu (\Delta h)^2
\ee
If we let $Y=\int_X|\Delta h|^2\o^n=\int_X|\nabla\bar\nabla h|^2\o^n$, it follows that
\be
\dot Y
=
-2\int_X|\nabla\Delta h|^2\o^n
+2\int_X|\nabla\bar\nabla h|^2\Delta h\o^n
+\int_X(\Delta h)^2(-R+\mu(n+2))\o^n
\leq
C\,Y
\ee 
since $\Delta h=R-\mu n$ is bounded in absolute value by assumption. As we saw earlier, this differential inequality
together with the integrability of $Y(t)$ over $[0,\infty)$
implies that $Y(t)\to 0$ as $t\to\infty$. Substituting this 
back in the above equation, we see also that $\int_0^\infty\int_X|\nabla (\Delta h)|^2\o^n<\infty$.

\medskip
$\bullet$ It would be interesting to find additional conditions, such as the $k$-nonnegativity of the traceless curvature operator, which would preserve the non-negativity of the Ricci curvature under the K\"ahler-Ricci flow in higher dimensions.
It is already known that the positivity of the Ricci curvature is not preserved by the Ricci flow on complete manifolds of dimensions 4 or higher, thanks to the counterexamples constructed by L. Ni \cite{N}.

\end{document}